\documentclass[12pt]{amsart}
\usepackage{amsmath,amsfonts,amsthm,amscd,textcomp,times}
\newtheorem{theorem}{Theorem}
\newtheorem{conjecture}{Conjecture}
\newtheorem{claim}{Claim}
\newtheorem{proposition}{Proposition}
\newtheorem{lemma}{Lemma}
\newtheorem{definition}{Definition}
\newtheorem{corollary}{Corollary}
\newtheorem{remark}{Remark}
\numberwithin{equation}{section}
\numberwithin{theorem}{section}
\numberwithin{proposition}{section}
\numberwithin{lemma}{section}
\numberwithin{claim}{section}
\numberwithin{corollary}{section}
\usepackage{diagrams}
\newcommand{\bull}{\ensuremath{{}\bullet{}}}
 
\newcommand{\cpn}{\ensuremath{\mathbb{P}^{N}}}
\newcommand{\slnc}{\ensuremath{SL(N+1,\mathbb{C})}}
\newcommand{\dlb}{\ensuremath{\overline{\partial}}}
\newcommand{\dl}{\ensuremath{\partial}}
\newcommand{\ra}{\ensuremath{\longrightarrow}}
\newcommand{\om}{\ensuremath{\omega}}
\newcommand{\vp}{\ensuremath{\varphi}}
\newcommand{\vps}{\ensuremath{\varphi_{\sigma}}}
\newcommand{\oms}{\ensuremath{\omega_{\sigma}}}

\newcommand{\omd}{\ensuremath{\widehat{\om}}}

\newcommand{\vplt}{\ensuremath{\varphi_{\lambda(t)}}}
\newcommand{\cpnd}{\ensuremath{\widehat{\mathbb{P}^{N}}}}

\newcommand{\cn}{\ensuremath{\mathbb{C}^{N+1}}}
 
\begin{document}
\title[ Duality and Energy Asymptotics ]{Projective Duality and K-Energy Asymptotics }
\author{Sean Timothy Paul}
\email{stpaul@math.wisc.edu}
\address{Mathematics Department at the University of Wisconsin, Madison}
\subjclass[2000]{53C55}
\keywords{Discriminants, Resultants, K-Energy maps, Projective duality .}
\date{July 2, 2008}
 \vspace{-5mm}
\begin{abstract}{ Let  $X\hookrightarrow \cpn$ be a smooth, linearly normal $n$ dimensional subvariety. Assume that the projective dual of $X$ has codimension one with defining polynomial $\Delta_X$. In this paper the log of the norm of 
$\sigma\cdot\Delta_X$ is expressed as the restriction  to the Bergman metrics of an energy functional on $X$.  We show how, for smooth plane curves, this energy functional reduces to the standard action functionals of K\"ahler geometry.}
\end{abstract}
\maketitle
\tableofcontents 
\section{Introduction}
 Let $(X, \om)$ be an $n$ dimensional compact K\"ahler manifold . A fundamental problem in K\"ahler geometry is to decide when the Mabuchi energy (see section 3, (\ref{directformula}) for the definition ) , denoted by $\nu_{\om}$, is bounded from below on the space of K\"ahler potentials. Lower bounds on $\nu_{\om}$ are closely related to the existence of canonical metrics (e.g. K\"ahler Einstein, constant scalar curvature, and extremal metrics ) in the class $[\om]$. When $X=\mathbb{P}^1$ to say that $\nu_{\om}\geq 0$
  is equivalent to the famous Moser-Trudinger inequality which plays a decisive role in the Nirenberg problem of prescribing Gauss curvature on $S^2$. Fundamental contributions to the study of lower bounds on $\nu_{\om}$ and the existence and uniqueness of canonical metrics are due to Bando and Mabuchi, Tian, S.K. Donaldson, and X.X. Chen and Tian .

  This paper bears on the special situation when $[\om]$ is a \emph{Hodge class}. That is, we assume that $\om=c_1(L,h)$ where $L$ is an ample line bundle on $X$, and $h$ is a Hermitian metric on $L$. We may as well assume (by raising $L$ to a power which does not concern us ) that $X$ is a subvariety of a \emph{fixed} $\cpn$  and that $\om=\om_{FS}|_X$. $X$ needs to be embedded in $\cpn$ in a sufficiently general manner, for example it should not lie in a hyperplane, nor should it arise as a non trivial projection from some larger projective space. What takes care of this is \emph{linear normality}, which means that the restriction map
  \[(\mathbb{C}^{N+1})^{\vee} \ra H^0(X,\ \mathcal{O}_X(1)) \]
  is an \emph{isomorphism} .
  Throughout this paper we shall assume that $X$ has this property.
  
   Whenever $X\hookrightarrow \cpn$ we may map $G:=\slnc$ into the space of K\"ahler potentials by pulling back the Fubini-Study form\footnote{Throughout this paper $G$ denotes $\slnc$, the special linear group over $\mathbb{C}$.}.  Therefore we may restrict the Mabuchi energy to the image of $G$. A rough formulation of the problem which motivates the present work is as follows.
  \ \\
  \begin{center} \emph{Relate the geometry of the embedding $X\hookrightarrow \cpn$ to the restriction of $\nu_{\om}$ to $G$ . }\end{center}
\ \\
In order to clarify this problem we first consider the situation for \emph{hypersurfaces}. This is due to Tian (see \cite{kenhyp}). This result is our model for higher codimensions.\\
\ \\
\noindent{\textbf{Theorem}} (Tian \cite{kenhyp}) \emph{Let $X_F$ be a smooth\footnote{Tian actually considers hypersurfaces with at worst normal singularities .} hypersurface in $\mathbb{C}P^{n+1}$ with defining polynomial $F$. Assume that $d=\mbox{deg}(F)\geq 2$. Then for all $\sigma\in G$ we have}
\begin{align}\label{hyp}
d \nu_{\om}(\vps)=\frac{(n+2)(d-1)}{n+1}\log\left(\frac{||\sigma\cdot F||^2}{||F||^2}\right)+\Psi_B([\sigma\cdot F])\ . 
\end{align}
\emph{ The ``singular" term $\Psi_B$ is given by}
\begin{align}
\Psi_{B}([f]):= \int_{X_f}\log\left(\frac{\sum_{i=0}^{n+1}|\frac{\dl f}{\dl z_i}(z)|^2}{||f||^2||z||^{2(d-1)}}\right)\om^n_{FS} \quad \mbox{\emph{for all $f\in B\setminus \Delta$}}\ .
\end{align}
 \emph{Moreover, $\Psi_{B}$ extends to the locus of reduced hypersurfaces, is bounded above, and}\newline\emph{ $\lim_{f_i\ra f_{\infty}}\Psi_{B}([f_i])=-\infty$ if and only if $f_{\infty}$ is non-reduced; $B:=\mathbb{P}(H^0(\mathbb{C}P^{n+1},\mathcal{O}(d)))$, and $\Delta$ denotes the discriminant locus.}\\
\begin{remark}
The issue, therefore, is to understand in terms of projective geometry, the ``singular" term $\Psi_B$.
\end{remark}
Next, we consider the Aubin energy restricted to $G$. In this case there is a completely satisfactory answer provided by Zhang and independently by the author (see \cite{gacms} and \cite{zhang}). \\
\begin{definition}(The Cayley-Chow Form)\\
 Let $X \subset \cpn$ be an $n$ dimensional irreducible subvariety of $\cpn$ with degree $d$.  Then the Cayley-Chow form of $X$ is given by
\begin{align*} 
&Z_{X}:= \{L \in \mathbb{G}| L\cap X \neq \emptyset\}. \\
\ \\
&\mathbb{G}:= \mathbb{G}(N-n-1,\mathbb{P}^N) \ .
\end{align*}
\end{definition}
\begin{remark}
It is easy to see that $Z_{X}$ is an {irreducible} hypersurface of degree $d$ in $\mathbb{G}$.
Since the homogeneous coordinate ring of the Grassmannian is a UFD, any codimension one subvariety is given by the vanishing of a section  $R_{X}$ of the homogeneous coordinate ring\footnote{See \cite{tableaux}  pg. 140 exercise 7.}
 \begin{align*}
 \{\ R_{X}=0\ \}= Z_{X}  \ ; \  R_{X} \in \mathbb{P}H^{0}(\mathbb{G},\mathcal{O}(d)).
 \end{align*}
 Following the terminology of Gelfand, Kapranov, and Zelevinsky \cite{gkz} we call $R_{X}$ the \textbf{$X$-resultant }.
 \end{remark} 

\textbf{Theorem .} 
\emph{Let $X$ be a linearly normal $n$ dimensional subvariety of \cpn. Let $R_X$ denote the $X$-resultant}
\emph{ Then the Aubin energy restricted to $G$ is given as follows}
\begin{align}\label{cayleyenergy}
-(n+1)\deg(X)F^0_{\om}(\vps)=\log\left(\frac{||\sigma\cdot R_X||^2}{||R_X||^2}\right) \ .
\end{align}
\emph{$|| \ ||$ is an appropriate norm on the vector space of polynomials of degree deg($X$) on $\mathbb{C}^{(n+1)(N+1)}$.}

In view of (\ref{hyp}) and (\ref{cayleyenergy}) I propose the following ``working conjecture". For further discussion of this see section 8.

\begin{conjecture}
Let $X^n$ be a linearly normal smooth subvariety of $\cpn$ of degree $d\geq2$. Then there are nontrivial finite dimensional complex $G$ representations $E_{1}$ and $E_{2}$ together with nonzero vectors $v_j=v_j(X)\in E_j$ such that for all $\sigma \in G$ the following identity holds
\begin{align}\label{tworeps}
\begin{split}
&\nu_{\om}(\vps)=\kappa_1\log\left(\frac{||\sigma\cdot v_1||_{E_1}^2}{||v_1||_{E_1}^2}\right)-\kappa_2 \log\left(\frac{||\sigma\cdot v_2||_{E_2}^2}{||v_2||_{E_2}^2}\right)\ .\\
\ \\
&\kappa_j\in \mathbb{Q}_+ \ \mbox{for}\ j=1,2\ \mbox{and $|| \ ||_{E_j}$ is a norm on  $E_j$ \ .}\\
\end{split}
\end{align}
 \end{conjecture}
 \begin{remark}
 The notation $v(X)$ is meant to suggest that $X$ is ``encoded" by  $v$.
  \end{remark}
What should the $v_j$ be? Certainly one of them should be an $X$-resultant since this reflects the presence of $J_\om$ and $I_\om$ in the Mabuchi energy. The main question then is what is the projective counterpart of the log term (see (\ref{directformula}) ) . First,
it is important to observe that the formation of the $X$-resultant is valid for an algebraic \emph{cycle}. In particular, \emph{smoothness is not required to make sense out of $R_X$}. Corresponding to this, the Aubin energy does not involve any of the \emph{curvature} of $\om_{FS}|_X$.  In order to incorporate the curvature of $X$ we require an analog of $R_X$ which \emph{is} sensitive to the smoothness of $X$. That is, to the \emph{existence of} $T^{1,0}_X$. Fortunately there is such an object, namely the \emph{projective dual to $X$}. Recall that the dual of a smooth projective variety consists of the \emph{tangent hyperplanes} to the variety (see \cite{gkz}) .
We shall denote the dual of $X$ by $X^{\vee}$ (or sometimes by $\widehat{X}$). In general, $X^{\vee}$ is a \emph{hypersurface} in the dual projective space.  When $X$ has dual defect equal to zero, i.e. when $X^{\vee}$ has codimension one, the defining polynomial is denoted by $\Delta_X$. Following Gelfand, Kapranov, and Zelevinsky we call $\Delta_X$ the \emph{$X$-discriminant} . Inspired by (\ref{cayleyenergy}) we consider
\begin{align}\log\left(\frac{||\sigma\cdot \Delta_X||^2}{||\Delta_X||^2}\right) \ .\end{align}
Where $|| \ ||$ is the standard norm on polynomials.

Now the problem is to find the analog of the left hand side of (\ref{cayleyenergy}). In other words,

\begin{center}\emph{We must associate to $\Delta_X$ an energy functional on $X$.}\end{center}
 
 The solution to this problem is given by Theorem (\ref{main}).
In a subsequent article we shall address the issue of relating this energy to the Mabuchi energy.  

  Our approach to this problem is a synthesis of several ideas, due to Simon Donaldson, Gang Tian, George Kempf, and Jean Michel Bismut, Henri Gillet and Christophe Soul\'e . In \cite{infdet} Donaldson introduced the fundamental idea of using the determinant of the Dirac operator as a functional, and then subjected this functional to the calculus of variations. The variation of this functional is expressed in terms of Bott-Chern secondary classes.
 In \cite{bottchrnfrms} Tian attached an energy functional, the Donaldson functional,  to any \emph{formal linear combination} of Hermitian holomorphic vector bundles on $X$ (see (1.9) pg. 212 of \cite{bottchrnfrms} ).   
  Tian then realized the Mabuchi energy as such a functional, and so introduced the notions of K stability and CM stability into K\"ahler geometry.  Bismut, Gillet and Soul\'e exhibited the Quillen norm of the canonical section of the determinant of the derived direct image of an acyclic \emph{complex} $(\xi^{\bull},\ \delta^{\bull})$ of holomorphic Hermitian vector bundles as the integral over the fiber of the double transgression of the Chern character of the complex multiplied by the Todd class of the fiber (see Theorem 0.3 pg. 51 of  \cite{bgs1}). In his study of collapsing, George Kempf introduced the \emph{Geometric Technique} for the study of syzygies (see section 4 of \cite{kempf76} and \cite{kempf73}). This technique consists in analyzing the derived direct image of the structure sheaf of a desingularization of the variety under study. In our case this variety is $X^{\vee}$. 
  
In a previous paper the author gave a new proof of the following result of Gelfand, Kapranov, and Zelevinsky \footnote{However, one should see Cayley's astonishing note on resultants 
 (\cite{cay}) and the work of Grothendieck, Knudsen and Mumford (see \cite{detdiv}) . } . 
 \begin{theorem}\label{disc&res} Let $X$ be a linearly normal smooth\footnote{Smoothness is only required for part $b)$. For $a)$ it is enough that $X$ be irreducible.} subvariety of $\cpn$. Let $\mathcal{V}$ be a holomorphic vector bundle over $X$. Let $({E}^{\bull}_{R}(\mathcal{V}(m)),\ \dl^{\bull}_{f})$ and $({E}^{\bull}_{\Delta}(\mathcal{V}(m)),\ \dl^{\bull}_{f})$ denote the resultant complex and the discriminant complex twisted by $\mathcal{V}(m)$ respectively. Then the following holds, provided $m\in \mathbb{Z}$ is sufficiently positive. 
 \begin{align*}
&a) \ \mbox{{The determinant of the resultant complex is the $X$-resultant .}}\\
&\qquad \mbox{\textbf{{Tor}}}({E}^{\bull}_{R}(\mathcal{V}(m)),\ \dl^{\bull}_{f})= R^{rank(\mathcal{V})}_X(f)\ ,\quad f\in  M_{n+1, N+1}(\mathbb{C})\ .\\
\ \\
&b) \ \mbox{Assume that the dual of $X$ is non-degenerate. Then the determinant}\\
&\mbox{of the discriminant complex is the $X$-discriminant .}\\
&\qquad \mbox{\textbf{{Tor}}}({E}^{\bull}_{\Delta}(\mathcal{V}(m)),\ \dl^{\bull}_{f})= \Delta^{rank(\mathcal{V})}_X(f)\ ,\quad  f\in (\mathbb{C}^{N+1})^{\vee}\ .
\end{align*}
\end{theorem}

Theorem \ref{disc&res} is an instance of the {Geometric Technique} . For a complete account the reader should consult \cite{gkz} and \cite{weyman}.
The fundamental ideas are due to A.Cayley (see \cite{cay}), A.Grothendieck\footnote{The details of Grothendieck's approach were written up in \cite{detdiv} and SGA 6 .}, G.Kempf, and I.M. Gelfand, M. Kapranov and A. Zelevinsky.
As this method is so vital for our present application we state it below. The precise aspects of the method that we require will be discussed in section 4 . \newpage 
  \begin{center}\emph{The Geometric Technique for Resultants and Discriminants}\end{center}
 
\begin{enumerate}
 \item Choose $a,b\ \in \mathbb{N}$ satisfying $(a+1)(N-b)=n+1$ . There are two cases. For $X$-resultants choose $a=0$ and $b=N-n-1$. For $X$-discriminants choose $a=n$ and $b=N-1$.
\ \\
\item Associate to $X$ an irreducible algebraic hypersurface $Z_{a,b}$ of an affine space $\mathbb{C}^k$.
$k(0,N-n-1)=(n+1)(N+1), \ k(n,N-1)= N+1$ . $Z_{0,\ N-n-1}=\{R_X=0\}$ and $Z_{n,\ N-1}=\{\Delta_X=0\}$. $Z_{a,b}$ is in general singular.
\ \\
\item Construct $I \hookrightarrow X\times \mathbb{C}^k$ a desingularization of $Z_{a,b}$ as an incidence correspondence. $I$ has the structure of a subbundle $\mathcal{S}$ of the trivial bundle $\mathcal{E}:=X\times \mathbb{C}^k$. 
\ \\
\item Let $\mathcal{Q}=\mathcal{E}/\mathcal{S}$. Form the Cayley Koszul resolution $K^{\bull}_{a,b}=K^{\bull}_{a,b}(\mathcal{Q})\rightarrow \iota_*(\mathcal{O}_I)\rightarrow 0$ of the structure sheaf of the desingularization $I$.
\ \\
\item Let $\mathcal{V}$ be any holomorphic vector bundle on $X$. Choose $m\in\mathbb{Z},\ m>>0$ . Compute $\mbox{\textbf{det}}R_{p_*}^{\bull}(K^{\bull}_{a,b}(\mathcal{V}(m)))$ .   
\end{enumerate}
\begin{remark}
{Theorem \ref{disc&res} states that $\mbox{\emph{\textbf{det}}}R_{p_*}^{\bull}(K^{\bull}_{a,b}(\mathcal{V}(m)))$  coincides with a power of the defining polynomial for $Z_{a,b}$} and that this power is equal to the rank of $\mathcal{V}$ .
\end{remark}

The basis of this paper consists in replacing step (5) of the Geometric Technique with \\
 
 \begin{center}(5)$^*$ \textbf{\emph{Extend the complex  $K^{\bull}_{a,b}$ to $G\times X$ and compute}} ${{p}_*}{\widetilde{\mbox{Ch}}}^{\bull} (K^{\bull}_{a,b}(\mathcal{V}(m)))$ . \end{center}
  \ \\
Here ${\widetilde{\mbox{Ch}}}^{\bull} (K^{\bull}_{a,b}(\mathcal{V}(m)))$  denotes the double transgression\footnote{This notion is explained in section 5 .} of the Chern character of the complex $K^{\bull}$, $p$ denotes the projection onto $G$, and $p_*$ denotes integration over the fiber. \\
\begin{remark}
This technique originated with George Kempf in his study of rational singularities of affine varieties associated to representations of  algebraic groups (see \cite{kempf76}). The method can be modified to compute the syzygies of the generic determinantal variety, symmetric and skew symmetric matrices, and nilpotent orbit closures . This is explained in \cite{weyman}.
\end{remark}
\begin{remark}{Theorem (\ref{main}) below says that the direct image of the double transgression of the Chern character of the complex $K^{\bull}_{a,b}(\mathcal{V}(m))$ gives the logarithm of the norm squared of the defining polynomial of} $Z_{a,b}$ . \end{remark}

\section{Statement of Results}
The main result of this paper is the following. 
  \begin{theorem}\label{main}
 Let $X\hookrightarrow \cpn$ be a smooth, linearly normal $n$ dimensional subvariety. Let $X^{\vee}$ be the dual of $X$. Assume that $X^{\vee}$ is a hypersurface with defining polynomial  $\Delta_X$.  Then there is a complex $(K^{\bull}, \dl^{\bull})$ of holomorphic vector bundles on $X$ and a smooth map
 $H:G\rightarrow \mathcal{M}_{K^{\bull}}$ such that
 \begin{align}\label{basicformula}
 D_{K^{\bull}}(\mbox{Ch}_{n+1}; H^{\bull}(e), H^{\bull}(\sigma))=\log\frac{||\sigma \cdot \Delta_X||^2}{||\Delta_X||^2} \ .
 \end{align}
 $D_{K^{\bull}}(\mbox{Ch}_{n+1}; H^{\bull}(e), H^{\bull}(\sigma))$ denotes the Donaldson functional \footnote{Precise definitions of all the terms appearing in the statement are given in subsequent sections.} of the complex $K^{\bull}$ .  $\sigma \in G$ ($e$ denotes the identity in $G$) , $|| \ ||$ is a norm on the vector space of degree $d^{\vee}:=\deg(X^{\vee})$ polynomials on $(\mathbb{C}^{N+1})^{\vee}$ . $\mathcal{M}_{K^{\bull}}$ denotes the space of smooth Hermitian metrics on $K^{\bull}$, and $Ch$ denotes the Chern Character .  
 \end{theorem}

This result exhibits the right hand side of (\ref{basicformula}) as the restriction of an energy functional on $X$ to the Bergman metrics associated to the embedding $X\hookrightarrow \cpn$ .

 \begin{remark}
We point out to the reader that ${{p}_*}{\widetilde{\mbox{Ch}}}^{\bull}(K^{\bull})=D_{K^{\bull}}$. The difference in notation reflects a difference in perspective. The notation $D_{K^{\bull}}$ is meant to emphasize the fact that the left hand side of (\ref{main}) is a functional on the space of K\"ahler potentials as in ( \cite{bottchrnfrms}). Whereas the double transgression notation emphasizes that we are taking the Quillen \footnote{See Theorem 0.3 pg. 51 of \cite{bgs1} .}norm of the canonical section of the determinant of the direct image of the complex $K^{\bull}$.
\end{remark}
\begin{corollary}
$\Delta_X$ is stable if and only if $G\ni \sigma \ra D_{K^{\bull}}(\mbox{Ch}_{n+1}; H^{\bull}(e), H^{\bull}(\sigma))$
is proper.
\end{corollary}

The next proposition identifies the energy  $D_{K^{\bull}}(\mbox{Ch}_{n+1}; H^{\bull}(e), H^{\bull}(\sigma))$ for smooth plane curves. Recall that
\[E_{1,\om}(\vp):= \int_X\log\left(\frac{{\omega_{\varphi}}}{\omega}\right)\left(\mbox{Ric}({\om})+\mbox{Ric}(\om_{\vp})\right)\]
is the \emph{Liouville energy} of the curve, and $F^0_{\om}$ is the Aubin energy\footnote{We recall the definitions of the basic energies in the next section . } 
\[F^0_{\om}(\vp):= J_{\om}(\vp)-\frac{1}{V}\int_X\vp\ \om^n \ . \] 
\begin{proposition}
Let $X_F:=\{F=0\} \hookrightarrow \mathbb{P}^2$ be a smooth nonlinear plane curve. Let $\Delta_F$ denote the discriminant of $F$. Then the energy functional is given by the formula
\begin{align} 
\begin{split}
D_{K^{\bull}}(Ch_2,H^{\bull}(e),H^{\bull}(\sigma))&=  4\deg(F)\nu_{\om}(\varphi_{\sigma}) -deg(F)E_{1,\om}(\varphi_{\sigma})\\
\ \\
&-4\deg(\Delta_F)F_{\om}^0(\vps) \ .
\end{split}
\end{align}
\end{proposition}
\begin{corollary}
 For all $\sigma \in SL(3, \mathbb{C})$ we have
\begin{align}\label{planecurve}
\begin{split}
  4\deg(F)\nu_{\om}(\varphi_{\sigma})&= \log \left(\frac{{||\sigma\cdot\Delta_{F} ||}^{2}}{{||\Delta_{F}||}^{2}}\right)- 2\frac{\deg(\Delta_F)}{\deg(F)}\log \left(\frac{{||\sigma\cdot{F}||}^{2}}{||{F}||^2}\right)\\
 \ \\ 
 &+\deg(F)E_{1,\om}(\varphi_{\sigma}) \ . 
\end{split}
\end{align}
\end{corollary}
 \section{Resume of Results of Part II}
 In a sequel to this paper we  verify our working conjecture for curves.
\begin{theorem}
 Let $X \hookrightarrow \cpn$ be a smooth, linearly normal curve of degree $d\geq 2$ . 
 Let   $R_{X}$  denote the \textbf{X-resultant} (the Cayley Chow form of $X$) . Let $\Delta_X$ denote the \textbf{$X$-discriminant} (the defining polynomial for the dual curve).
 Then the K-energy map restricted to the Bergman metrics is given as follows
\begin{align} \label{bergm}
\begin{split}
&d\nu_{\om}(\varphi_{\sigma})= \log \left(\frac{{||\sigma\cdot\Delta_{X}||}^{2}}{{||\Delta_{X}||}^{2}}\right)- \frac{\deg(\Delta_X)}{\deg(R_X)}\log \left(\frac{{||\sigma\cdot R_{X}||}^{2}}{||R_{X}||^2}\right)  \ .
 \end{split}
\end{align}
\end{theorem} 
 
( \ref{bergm}) is compatible with (\ref{planecurve}). In the sequel we may deduce the following corollary.

\begin{corollary}(Stability and the Moser-Trudinger inequality)\\
Let $X^1_2$  denote the second Veronese image of $\mathbb{P}^1$ . This is a quadric $\{Q=0\}$ in $\mathbb{P}^2$.   Let $\Delta_Q$ denote the dual quadric. Then for all $\sigma \in SL(3,\mathbb{C})$ we have the identity
\begin{align}
\log \left(\frac{{||\sigma\cdot\Delta_{Q} ||}^{2}}{{||\Delta_{Q}||}^{2}}\right)=6E_{1,\om}(\varphi_{\sigma}) \ .
\end{align}
Since $(X^1_2)^{\vee}$ is a smooth curve, and smooth hypersurfaces are automatically semistable, we have the inequality
\begin{align}
E_{1,\om}(\varphi_{\sigma})\geq C \quad \mbox{for all $\sigma\in SL(3,\mathbb{C})$}\ .
\end{align}
\end{corollary}
 \begin{remark}
 The point is that I do not need the Moser-Trudinger inequality in this case.  
 \end{remark}
 \begin{definition} (Discriminant and Chow Polytopes)
Let $X\hookrightarrow \cpn$ be a smooth linearly normal dually non-degenerate subvariety .Then the \textbf{discriminant polytope} of $X$ is the weight polytope $N(\Delta_X)$ of the $X$-discriminant.
The \textbf{Chow Polytope} is the weight polytope $N(R_X)$ of the $X$-resultant (the Cayley-Chow form of $X$).
\end{definition}
For the precise definition of the weight polytope, please see section 8. Below $l_{\lambda}$ denotes the linear functional on $\mathbb{R}^N$ corresponding to $\lambda$ , where $\lambda$ is an algebraic one parametrer subgroup of $G$.
 \begin{corollary}(Mabuchi Energy Asymptotics on Algebraic Curves)\\
Let $X\hookrightarrow \cpn$ be a smooth, linearly normal algebraic curve of degree $d\geq 2$. Then there is an asymptotic expansion as $|t|\rightarrow 0$
\begin{align}
d\nu_{\om}(\vplt)=\left(\mbox{\emph{Min}}_{ \{ x\in N(\Delta_X)\}}\ l_{\lambda}(x)- 
 \frac{\deg(\Delta_X)}{\deg(R_X)}   \mbox{\emph{Min}}_{ \{ x\in N(R_X)\}}\ l_{\lambda}(x)\right)\log(|t|^2)+O(1)\ .
\end{align}
\end{corollary}
The previous proposition is particularly interesting when $X$ is the \emph{rational normal curve} $X^1_d$ ($d\geq 2$). We state this explicitly in the following corollary. 
\begin{corollary} 
The Mabuchi energy is bounded below along all degenerations $\lambda$ if and only if
\begin{align}\label{polymoser}
 \left(\frac{d-1}{d}\right)N(R_{X_d^1})\subset N(\Delta_{X_d^1})\ .
\end{align}
\end{corollary}

 This paper is organized as follows. In section 3 we recall basic definitions and notations from K\"ahler geometry and give a brief account of the Bergman metrics. In section 4 we describe that part of the geometric technique most relevant to the construction of discriminants and resultants. In particular, the complex $K^{\bull}$ is defined. In section 5 we discuss the double transgression of the Chern character of a complex of holomorphic Hermitian vector bundles on a complex manifold. The proof of the main theorem takes up section 6. The argument consists in equipping the terms of the complex $K^{\bull}$ with Hermitian metrics depending on $G$ and applying the double transgression construction of section 5. In section 7 we consider the example of smooth plane curves. In this case $D_{K^{\bull}}$ is a linear combination of the Mabuchi energy, the Aubin energy, and the Liouville energy. Section 8 relates the author's working conjecture to the asymptotic behavior of the Mabuchi energy along arbitrary degenerations.  
 
\section{Preliminaries and Notations from K\"ahler Geometry}
Let $(X,\ \om)$ be a K\"ahler manifold. We always set $\mu$ to be the average of the scalar curvature of $\om$ and $V$ to be the volume
\begin{align*}
&\mu:=\frac{1}{V}\int_X\mbox{Scal}(\om)\om^n \\
\ \\
&V:= \int_X\om^n \ .
\end{align*}

The space of K\"ahler potentials will be denoted by $\mathcal{H}_\om$
\begin{align*}
\mathcal{H}_\om:=\{\vp\in C^{\infty}(X)|\om_\vp:=\om+\frac{\sqrt{-1}}{2\pi}\dl\dlb\vp>0 \} \ .
\end{align*}

The \emph{Mabuchi K-Energy}, denoted by $\nu_\om$,  is a map $\nu_\om:\mathcal{H}_\om\ra \mathbb{R}$ and is given by the following expression
\begin{align*}
 \qquad \nu_{\omega}(\varphi):= -\frac{1}{V}\int_{0}^{1}\int_{X}\dot{\varphi_{t}}(\mbox{Scal}(\varphi_{t})-\mu)\omega_{t}^{n}dt.
\end{align*}
Above, $\varphi_{t}$ is a smooth path in $\mathcal{H}_\om$ joining $0$ with $\varphi$. The K-Energy does not depend on the path chosen. If $\vp$ is a critical point of the Mabuchi energy then $\mbox{Scal}_\om(\vp)\equiv \mu$ .

Suppose that $\om$ satisfies $\mbox{Ric}(\om)= \frac{\mu}{n}\om+\frac{\sqrt{-1}}{2\pi}\dl\dlb h_\om$ .   In this case there is the following well known direct formula for the K-energy map. 

\begin{align} \label{directformula}
\begin{split}
&\nu_{\omega}(\varphi)=\int_{X}\mbox{log}\left(\frac{{\omega_{\varphi}}^{n}}{\omega^{n}}\right)\frac{{\omega_{\varphi}}^{n}}{V} - \frac{\mu}{n}(I_{\omega}(\varphi)-J_{\omega}(\varphi)) -\frac{1}{V}\int_Xh_\om(\om^n_\vp - \om^n) \\
\ \\
& J_{\omega}(\varphi):= \frac{1}{V}\int_{X}\sum_{i=0}^{n-1}\frac{\sqrt{-1}}{2\pi}\frac{i+1}{n+1}\dl\varphi \wedge \dlb
\varphi\wedge \omega^{i}\wedge {\omega_{\varphi} }^{n-i-1}\\
\ \\
&I_{\omega}(\varphi):= \frac{1}{V}\int_{X}\varphi(\omega^{n}-{\omega_{\varphi}}^{n})\ .
\end{split}
\end{align}

Assume that $\om=-\frac{\sqrt{-1}}{2\pi}\dl\dlb \log(h)$ where $L$ is a very ample line bundle on $X$ with Hermitian metric $h$. When $X$ is Fano, $L$ may be taken to be some large multiple of $-K_X$ . By definition, there is an embedding
\[
\begin{CD}
X@>>{ L^{}}>\mathbb{P}(H^{0}(X, L^{})^{*})= \mathbb{P}^{N} \ .
\end{CD}
\]
furnished by some basis $\{S_{0},\dots,S_{N}\}$ of $H^{0}(X, L)$.
In this paper the main concern is with a \emph{fixed} embedding, in which case we may as well take $L=\mathcal{O}_{\mathbb{P}^N}(1)|_X$ and correspondingly $\om=\om_{FS}$ . This amounts to setting the $S_i$ to be the restriction of homogeneous coordinates to $X$ .

Let $\sigma \in G$, then
\begin{align}
\sigma^{*}(\omega_{FS})=\omega_{FS} +\frac{\sqrt{-1}}{2\pi}\dl\dlb\varphi_{\sigma} >0 \ .  
\end{align}
It is easy to see that  $\varphi_{\sigma}$ is given by the formula
\begin{align}
\varphi_{\sigma}=\mbox{log}\left(\frac{||\sigma z||^{2}}{|| z||^{2}}\right).
\end{align}

Then the \emph{Bergman metrics} associated to the embedded of $X$ in $\cpn$ are by definition 
\begin{align}
\mbox{Berg}_N:=\{\oms:=\om+\frac{\sqrt{-1}}{2\pi}\dl\dlb\vps|\ \sigma\in G\}\hookrightarrow \mathcal{H}_\om \ .
\end{align}

 Now we may \emph{consider the K-energy as a function on $G$} .  
 
 In the text $J_\om(\sigma)$, $I_\om(\sigma)$, and $\nu_\om(\sigma)$ denote $J_\om(\vps)$ , etc. .
  
\section{The complex $K^{\bull}$}
   Our concern is with irreducible subvarieties $Z$ of an affine space $\mathbb{C}^{{k}}$ associated to a smooth, linearly normal subvariety $X$ of $\cpn$. Such subvarieties $Z$ arise in the following manner.
\emph{Assume} there exists a vector subbundle $\mathcal{S}$ of the trivial bundle $\mathcal{E}:= X\times \mathbb{C}^{k}$ such that the image of the restriction to $I$  of the projection of $\mathcal{E}$ onto  $\mathbb{C}^{k}$ is $Z$, where $I$ denotes the \emph{total space} of $\mathcal{S}$. We shall always take $f$ to be a variable point in $\mathbb{C}^{k}$.
 
 There is the exact sequence of vector bundles on $X$
\begin{align*}
0\rightarrow \mathcal{S}\rightarrow \mathcal{E}\overset{\pi}\rightarrow \mathcal{Q}\rightarrow 0 \ .
 \end{align*}
In this case there is tautological \emph{regular} section $s$ of ${p_1}^*(\mathcal{Q})$ whose base locus is $I$.   $p_1$ denotes the projection of $\mathcal{E}$ to $X$.
We let $p_I$ denote the restriction of $p_2$ to $I$. $Z$ denotes the image of $I$ under $p_I$. This situation is pictured below in what we will call the \emph{ basic set up} following the terminology of J.Weyman (see \cite{weyman} ).
\begin{diagram}
&&{p_1}^*(\mathcal{Q})&\rTo^{\pi_2}&\mathcal{Q}\\
&& \dTo^{\pi_1}&&\dTo^{p}\\
 I&\rTo^{\iota}&{X}\times  \mathbb{C}^{k}&\rTo^{p_1}&X\\
\dTo ^{p_{I}}&&\dTo^{p_2}\\
  Z&\rTo ^{i}& \mathbb{C}^{k}&
\end{diagram}

\emph{In our applications we shall have that $Z$ is an irreducible algebraic {\textbf{hypersurface}} in $\mathbb{C}^{k}$, and that $p_I:I\rightarrow Z$ is a resolution of singularities. Therefore, in the remainder of this section we assume that $Z$ has codimension one .} 

Observe that the assumption on the codimension of $Z$ in $\mathbb{C}^k$ forces $\mbox{rank}(Q)=n+1$. 
In this case, following G. Kempf (see the section on ``Historical Remarks" in (\cite{kempf76}) ), we may study the irreducible equation of $Z$ through an analysis of the direct image of a Cayley-Koszul complex of sheaves on ${X}\times \mathbb{C}^{k}$. We have the free resolution over $\mathcal{O}_{X\times \mathbb{C}^{k}}$
\begin{align}
(K^{\bull}(p_1^*(\mathcal{Q}^{\vee})), (s\wedge \cdot)^*)\rightarrow \iota_{*}\mathcal{O}_{I}\rightarrow 0\ ; \ K^j(p_1^*(\mathcal{Q}^{\vee})):=\bigwedge^{n+1-j}p_1^*(\mathcal{Q}^{\vee})\ .
\end{align}
More generally, let $\mathcal{V}$ denote any vector bundle on $X$. Then we will consider the \emph{twisted} 
 complex
\begin{align}
\begin{split}
&(K^{\bull}(p_1^*(\mathcal{Q}^{\vee}))\otimes p_1^*\mathcal{V} , \ (s\wedge \cdot)^*)\rightarrow \iota_{*}\mathcal{O}_{I}\otimes  p_1^*\mathcal{V}\rightarrow 0 \\
\ \\
&(s\wedge \cdot)^* \ \mbox{denotes interior multiplication} \ .
\end{split}
\end{align}
 Let $f\in  \mathbb{C}^{k}$, then we may pull the twisted Cayley-Koszul complex back to $X$ via the map
\begin{align}
i_{f}:X\rightarrow X\times \mathbb{C}^{k} \quad i_f(x):= (x,f)
\end{align}

\begin{center}\emph{Then $i_f^*(K^{\bull}(p_1^*(\mathcal{Q}^{\vee}))\otimes p_1^*\mathcal{V} , \ (s\wedge \cdot)^*)$ is an 
 {acyclic}  complex of vector bundles on $X$ whenever  $f\in \mathbb{C}^{k}\setminus Z$.}\end{center}

Given $X\hookrightarrow \cpn$ we can achieve the situation of the basic set up as follows. Let $a, \ b \in \mathbb{N}$ satisfy $(a+1)(N-b)=n+1$. There is a map $\iota_X$ from $X$ into $\mathbb{G}(a,\ N)$ ( the Grassmannian of $a$ dimensional linear subspaces of $\cpn$ ) given by inclusion when $a=0$ or the Gauss map when $a=n$. In either case, this map is \emph{finite to one} as follows from celebrated work of F. L. Zak (see \cite{zak}).  
\begin{diagram}
&&   ({\iota_X\times\mathbb{I}_{\mathbb{G}}})^* (\mathcal{U}^{\vee}_{a+1}\otimes\mathcal{Q}_{N-b})&\rTo^{\pi_2}&\mathcal{U}^{\vee}_{a+1}\otimes\mathcal{Q}_{N-b}\\
&& \dTo^{\pi_1}&&\dTo^{p}\\
 I_{\mathbb{G}}&\rInto^{\iota}&{X}\times  \mathbb{G}(b,N)&\rTo^{\iota_X\times\mbox{id}}&\mathbb{G}(a,N)\times \mathbb{G}(b,N)\\
\dTo ^{p_{I}}&&\dTo^{p_2}\\
  Z_{\mathbb{G}}&\rInto ^{i}&\mathbb{G}(b,N)  &
\end{diagram}
 There is a natural section $s$ of the bundle $\mathcal{U}^{\vee}_{a+1}\otimes\mathcal{Q}_{N-b}$ over $\mathbb{G}(a,N)\times \mathbb{G}(b,N)$ given by
 \begin{align*}
 s|_{(E,L)}:E\subset \mathbb{C}^{N+1}\rightarrow \mathbb{C}^{N+1}/L \ .
\end{align*}
 $I_{\mathbb{G}}$ denotes the base locus of $ (\iota_X\times\mbox{id})^*(s)$. $p_I$ is the restriction of the second projection. In most cases, $Z:= p_I(I_{\mathbb{G}})$ has codimension one in $\mathbb{G}(b,N)$. 
 
 To make contact with the affine situation we proceed as follows.   
 The affine space is defined by $\mathbb{C}^{k}:= M_{(N-b)\times (N+1)}(\mathbb{C})$.
Let $r$ denote the rational projection map from $\mathbb{C}^{k}$ to $\mathbb{G}(b,N)$
 \begin{align*}
 \mathbb{C}^{k} \overset{r}{\dasharrow} \mathbb{G}(b, \ N)\ .   
  \end{align*}
The incidence correspondence $I$ is defined to be 
  \begin{align*}
  I:= \overline {(1_X\times r)^{-1}(I_{\mathbb{G}})}=\{(y,\ A)\in X\times \mathbb{C}^k|\iota_X(y)\subset \mbox{ker}(A) \}\ .
  \end{align*}
  $1_X\times r$ denotes the map 
\begin{align*}
X\times \mathbb{C}^{k}\overset{1_X\times r }{\dasharrow}X\times \mathbb{G}(b, \ N)\ .
\end{align*}
 $I$ has the structure of a subbundle of $\mathcal{S}$ of  the trivial bundle $X\times \mathbb{C}^{k}$, and we have that
 \begin{align*}
Q:=\mathcal{E}/\mathcal{S} \cong \iota_{X}^*(\mathcal{U}_{a+1}^{\vee})\otimes \mathbb{C}^{N-b} \ .
\end{align*}
  The image of the restriction to $I$ of the projection $p_2$ is given by 
  \begin{align*}
  Z:=\overline{r^{-1}(Z_{\mathbb{G}})}\ .
\end{align*}
\begin{definition} Let $f\in M_{(N-b)\times (N+1)}(\mathbb{C})$, be generic. Then we define an acyclic complex  $(K_{a,b}^{\bull},\ \dl_{f}^{\bull})$ of locally free sheaves on $X$ as follows.
 \begin{align}
 (K_{a,b}^{\bull}, \ \dl_{f}^{\bull}):= \left(\bigwedge^{n+1-\bull}(\iota_X^*\mathcal{U}_{a+1}\otimes \mathbb{C}^{N-b})(m)|_{X\times \{f\}}, \ (s(\iota_X(\cdot ),\ f )\wedge)^*\right)
\end{align}
When $a=0$ this complex is denoted by $(K_R^{\bull}(m),\ \dl_f^{\bull})$, when $a=n$, the case considered in this paper, the complex shall be denoted by $(K_{\Delta}^{\bull}(m),\ \dl_f^{\bull})$ . 
\end{definition}
We may write, 
\begin{align}
(K^{\bull}_{\Delta}(m), \ \dl_f^{\bull})= (\bigwedge^{n+1-\bull}(J_1(\mathcal{O}_X(1))^{\vee})(m)|{X\times \{f\}},\ (s(\rho_X(\cdot), f)\wedge \cdot)^*) 
\end{align}
$J_1(\mathcal{O}_X(1))$ denotes the bundle of \emph{one jets}. This a rank $n+1$ bundle over $X$.  We recall the definition. 

To begin we consider the \emph{affine cone} over $X$, which we denote by $\tilde{X}$.  $\tilde{X}$ is a smooth subvariety of $\mathbb{C}^{N+1}\setminus \{ 0\}$.
Let $\{F_{\alpha}\}$ denote any generating set for the homogeneous ideal of $X$. Then
\begin{align*}
J_1(\mathcal{O}_X(1))^{\vee}:=T^{1,0}(\tilde{X})=\{(p,w)\in X\times \mathbb{C}^{N+1}|\ \nabla F_{\alpha}(p)\cdot w=0 \ \mbox{for all $\alpha$} \}\overset{\iota}{\hookrightarrow} X\times \mathbb{C}^{N+1}\ .
\end{align*}
In section 7 we will require the following well known fact.
 \begin{proposition} \label{jetseq}There is an exact sequence of vector bundles on $X$.
 \begin{align}\label{xact}
 \begin{split}
 0\rightarrow \mathcal{O}_X(-1)\overset{\iota}{\rightarrow}&T^{1,0}(\tilde{X}) \overset{\pi}{\rightarrow} T^{1,0}(X)\otimes \mathcal{O}_X(-1)\rightarrow 0 
 \end{split}
 \end{align}
 \end{proposition}

\begin{remark}
 \emph{On the one hand $\pi$ denotes the map
\begin{align*}
T^{1,0}(\tilde{X})\overset{\pi}{\rightarrow} T^{1,0}(X)\otimes \mathcal{O}_X(-1)\rightarrow 0 \ .
\end{align*}
On the other hand we \emph{also} denote by $\pi$ the projection onto $\cpn$
\begin{align*}
\pi:\mathbb{C}^{N+1}\setminus \{0\}\rightarrow \cpn \ .
\end{align*}
Finally we can define $\pi$ in (\ref{xact}) by the formula (where $\pi(v)=p$ )}
\begin{align*}
T^{1,0}(\tilde{X})  \ni (p,w)\rightarrow \pi(p,w):= {\pi_{*}}|_v(w)\otimes v \in T^{1,0}(X)\otimes \mathcal{O}_X(-1)\ .
\end{align*}
\end{remark}
For fixed $p\in X$, $T_p^{1,0}(\tilde{X})\subset \cn $ is a linear subspace the projection of which is the \emph{embedded tangent space to $X$ at $p$} . The Gauss map $\rho_X$ is given by 
\begin{align}\label{embedded}
\rho_X(p):=\mathbb{P}(T_p^{1,0}(\tilde{X}))\subset \cpn \ .
\end{align}
Often we use the notation $\mathbb{T}_p(X)$ to denote the embedded tangent space to $X$ at the point $p$.
\section{Bott-Chern Forms}
Let $\phi$ be a $GL_n(\mathbb{C})$ invariant polynomial on $M_{n\times n}(\mathbb{C})$ homogeneous of degree $d$. The \emph{complete polarization} of $\phi$ is defined as follows. Let $\tau_1,\tau_2,\dots, \tau_d$ be arbitrary real parameters. Then
\begin{align*}
\phi(\tau_1A_1+\tau_2A_2+\dots +\tau_dA_d)=\sum_{|\alpha|=d}\phi_{\alpha}(A_1,A_2,\dots,A_d)\tau^{\alpha} \  \quad \tau^{\alpha}:=\tau_1^{\alpha_1}\tau_2^{\alpha_2}\dots\tau_d^{\alpha_d} \ .
\end{align*}
We let $\phi_{(1)}(A_1,A_2,\dots,A_d)$ denote the coefficient of $\tau_1\tau_2\dots\tau_d$ . We define
\begin{align*}
\phi_{(1)}(A;B):=\phi_{(1)}(A,\overbrace{B,B,\dots,B}^{d-1}) \ .
\end{align*}
 Let $M$ be an $n$ dimensional complex manifold, $E$ is a holomorphic vector bundle of rank $k$ over $M$. $H_0$ and $H_1$ are two Hermitian metrics on $E$.
Let $H_t$ be a smooth path joining $H_0$ and $H_1$ in $\mathcal{M}_{E}$ ( the space of Hermitian metrics on $E$) . Define $U_t:= (\frac{\dl}{\dl t}H_t)\cdot H^{-1}_t $ . $F_t:= \dlb\{(\dl H_t)H^{-1}_t\}$ is the curvature of $H_t$ (a purely imaginary (1,1) form) .  Now suppose that $\phi$ is a homogeneous invariant polynomial on $M_{k\times k}(\mathbb{C})$ of degree $d$. Then 
\begin{align*}
\phi_{(1)}(U_t;\ F_t)
\end{align*}
is a form of type $(d-1, d-1)$ .
The Bott Chern form is given as follows 
\begin{align*}
BC(E,\phi ; H_0,H_1):=-\frac{\deg(\phi)}{(n+1)!}\int_0^1\phi_{(1)}(U_t;\ \frac{\sqrt{-1}}{2\pi}F_t)\ dt \ .
\end{align*}
\begin{proposition} (Bott and Chern)
\begin{align*}
\frac{\sqrt{-1}}{2\pi} \dl\dlb BC(E,\phi ; H_0,H_1)=\phi(\frac{\sqrt{-1}}{2\pi}F_1)-\phi(\frac{\sqrt{-1}}{2\pi}F_0) \ .
\end{align*} 
\end{proposition}
When ${\deg}(\phi)$ has degree $n+1$ $BC(E,\phi ; H_0,H_1)$ is a \emph{top dimensional} form on $M$, and the following integral is well defined
\begin{align*}
D_{E}(\phi; H_0,H_1):= \int_M BC(E,\phi ; H_0,H_1) \ .
\end{align*}
Of particular importance is when $\phi (A)=Ch_{n+1}(A)=\frac{1}{(n+1)!}\mbox{Tr}(A^{n+1})$ . Observe that in this case we have
\begin{align*}
 \phi_{(1)}(A;B)=\mbox{Tr}(AB^n) \ .
\end{align*}
Let $H:Y\rightarrow \mathcal{M}_{{E}}$ (the space of $C^{\infty}$ Hermitian metrics on $E$)  be a smooth map, where $Y$ is a complex manifold of dimension $m$. \emph{Fix} a Hermitian metric $H_0$ on $E$.
Then we are interested in the smooth function on $Y$
\begin{align*}
Y\ni y \rightarrow D_{E}(\phi; H_0, H(y)) \ .
\end{align*}

Let $p_2$ denote the projection from $Y\times M$ onto $M$. Then $H(y)$ is a smooth Hermitian metric on $p^{*}_2(E)$  
 whose curvature is given by
\begin{align*}
F_{Y\times M}(H(y))=\dlb_{Y\times M}\{(\dl_{Y\times M}H(y))H(y)^{-1}\} \ .
\end{align*}
For the proof of the following proposition, see \cite{bottchrnfrms} prop. 1.4 on pg. 213 .
\begin{proposition} Let $\phi$ be homogeneous of degree $n+1$ and $H_0$ a fixed metric on ${E}$. Then for all smooth compactly supported forms $\eta$ of type $(m-1,m-1)$ we have the identity
\begin{align}\label{workhorse}
\frac{\sqrt{-1}}{2\pi}\int_YD_{E}(\phi; H_0, H(y))\dl_Y\dlb_Y\eta=\int_{Y\times M}\phi(F_{Y\times X}(\frac{\sqrt{-1}}{2\pi}H(y)))\wedge p_1^*(\eta)\ .
\end{align}
\end{proposition}

Now we extend the above to holomorphic Hermitian \emph{complexes} $(E^{\bull}, H_0^{\bull} ; \dl^{\bull})$ .
Let $H^{\bull}:Y\rightarrow \mathcal{M}_{E^{\bull}}$ be a smooth map. Concretely, $H^{\bull}(y)$ is a $C^{\infty}$ metric on $E^{\bull}$. We define the smooth function on $Y$.
\begin{align*}
D_{E^{\bull}}(\mbox{Ch}_{n+1}; H_0^{\bull}, H^{\bull}(y)):=
\sum_{j=0}^{l}(-1)^jD_{E^{j}}(\mbox{Ch}_{n+1}; H_0^{j}, H^{j}(y)) \ .
\end{align*}

\begin{corollary}
For all smooth compactly supported $(m-1,m-1)$ forms $\eta$ on $Y$ we have
\begin{align*}
&\int_Y\frac{\sqrt{-1}}{2\pi}D_{E^{\bull}}(\mbox{Ch}_{n+1}; H_0^{\bull}, H^{\bull}(y))\wedge \dl_Y\dlb_Y \eta=\\ 
\ \\
&\sum_{j=0}^{l}(-1)^j\int_{Y\times M}\mbox{Ch}_{n+1}(\frac{\sqrt{-1}}{2\pi}F_{Y\times X}^{E^j}(H^j(y))\wedge p_1^*(\eta) \ .
\end{align*}
\end{corollary}
 
 \section{Proof of the Main Theorem}
  The first step is to reinterpret the complexes $(K^{\bull}_{\Delta}(m),\ \dl^{\bull}_{f})$ (and  ($K^{\bull}_{R}(m),\ \dl^{\bull}_{f})$) as complexes of sheaves on $GX$ where
  \[GX:=\{(\sigma, y)\in G\times \cpn| \ y \in \sigma X \}\ .\]   
To carry this out for discriminants recall that the \emph{Gauss map} associated to $X\hookrightarrow \cpn$ is given by
\[ 
\rho: X \rightarrow \mathbb{G}(n,\mathbb{P}^{N})\qquad \rho (p)=\mathbb{T}_{p}(X) \ .
\]
$\mathbb{G}(n,\mathbb{P}^{N})$ is the Grassmannian of $n$ dimensional linear subspaces of $\mathbb{P}^{N}$ and $\mathbb{T}_{p}(X)$ denotes the \emph{embedded tangent space} to $X$ at $p$ (see \ref{embedded}) .
Let $f$ be a linear form on $X$, that is, $f \in H^{0}(X,\mathcal{O}_{X}(1))$ we define a map $\rho_{G,f}$ as follows.
\[
\rho_{G,f}: GX \rightarrow \mathbb{G}(n,\mathbb{P}^{N})\times {\mathbb{P}^{N}}^{\vee}\qquad \rho_{G,f}(\sigma,y):= (\mathbb{T}_{y}(\sigma X),\sigma f)
\]
 Let $\mathcal{U}$ denote the rank=$n+1$ universal vector bundle on $\mathbb{G}(n,\mathbb{P}^{N})$. Consider the bundle 
\begin{align} 
 F:= p_{1}^{*}\mathcal{U}^{\vee}\otimes p_{2}^{*}\mathcal{O}_{{\mathbb{P}^{N}}^{\vee}}(1)
 \end{align}
over the product $\mathbb{G}(n,\mathbb{P}^{N})\times {\mathbb{P}^{N}}^{\vee} $. 
There is a canonical regular section $s$ of this bundle whose base locus is a flag manifold
 \begin{align}\label{bslocus}
\begin{split}
&I := \{(L, f)\in \mathbb{G}(n,\mathbb{P}^{N})\times {\mathbb{P}^{N}}^{\vee}|\  L\subset \mbox{ker}(f) \}\overset{\iota}\rightarrow \mathbb{G}(n,\mathbb{P}^{N})\times {\mathbb{P}^{N}}^{\vee} \\
\ \\
&I =\{s=0\} \ .
\end{split}
 \end{align}
 Just as in the affine case  $\iota _{*}\mathcal{O}_{I_{\Delta}}$ is resolved by the Cayley-Koszul complex  
\begin{align*}
(\bigwedge^{n+1-\bull}F^{\vee},\ \dl^{\bull}:= (s\wedge \cdot)^* )\ .
\end{align*}
On $GX$ we introduce the following  complex, where $f\in (\mathbb{P}^N)^{\vee}$ is chosen generically .  
\begin{align}\label{basicmplx}
\begin{split}
&(K^{\bull}_{G\Delta}(m),\ \dl^{\bull}_{f}) := \left( \rho_{G,f}^{*}(\bigwedge^{n+1-\bull}F^{\vee})\otimes \pi^{*}\mathcal{O}_{X}(m),\   (s\circ \rho_{G,f}\wedge \cdot)^*  \right) \\
 \end{split}
\end{align}
 $\pi$ denotes the projection of $GX$ onto $X$.  
We note that under the composition
\begin{align*}
X\overset{\iota_e}{\hookrightarrow}GX\overset{\rho_{G,f}}{\rightarrow} \mathbb{G}\times {\mathbb{P}^N}^{\vee}
\quad {\iota_e}(x):= (e,x) \ 
\end{align*}
the complex $(\bigwedge^{n+1-i}F^{\vee},\dl^{\bull})$ pulls back to $({K}^{\bull}_{\Delta}(m),\ \dl^{\bull}_{f})$.
Recall that when $f\notin \widehat{X}$ the complex  $(K^{\bull}_{G\Delta}(m),\ \dl^{\bull}_{f})$ is exact. Moreover this complex carries a natural \emph{Hermitian metric} (on each term) induced by the natural metrics on $U$ and $\mathcal{O}_{{\mathbb{P}^{N}}^{\vee}}(1)$.

 Before we proceed to the proof of the main theorem, let us explain what is meant by a continuous metric (or norm) on $\mathcal{O}_{B}(-1)$,  where $B:=\mathbb{P}(H^0(\widehat{\cpn}, \mathcal{O}(\widehat{d})))$ and $\widehat{d}$ denotes the degree of $X^{\vee}$ .  \emph{Up to scaling} we have that $\Delta_X \in H^{0}({\mathbb{P}^{N}}^{\vee},\mathcal{O}(\widehat{d}))$.

In general we write linear form $f$ on $\cpn$  (i.e. a \emph{point} in the \emph{dual} $\cpn$) as $f= a_0z_0+a_1z_1+\dots +a_Nz_N$. Therefore we take $[a_0:a_1:\dots :a_N]$ as the homogeneous coordinates of $f$ on $\widehat{\cpn}$.
Therefore we may write
\begin{align*}
\Delta_X(f)=\sum_{|\alpha|=\widehat{d}}c_{\alpha_0,\dots \alpha_N}{a_0}^{\alpha_0}{a_1}^{\alpha_1}\dots {a_N}^{\alpha_N} \ .
\end{align*}
The \emph{finite dimensional complex vector space} $H^{0}({\mathbb{P}^{N}}^{\vee},\mathcal{O}(\widehat{d}))$ comes equipped  with its standard Hermitian inner product $< ,\ >$ in which the monomials ${a_0}^{\alpha_0}{a_1}^{\alpha_1}\dots {a_N}^{\alpha_N}$ form an orthogonal basis. Under a suitable normalization we have that
\begin{align*}
||\Delta_X||_{FS}^2:=<\Delta_X,\Delta_X>=\sum_{|\alpha|=\widehat{d}}\frac{|c_{\alpha_0,\dots \alpha_N}|^2}{\alpha_0!\alpha_1!\dots \alpha_N!}\ .
\end{align*}
Finally, to say that the metric $||\ ||$ on $\mathcal{O}_{B}(-1)$ is \emph{continuous} means that there is a continuous function $\theta$ on $B$ such that
\begin{align}
\exp({\theta})||\ ||_{FS}= ||\ || \ .
\end{align}
Since $B$ is compact, the conformal factor $\exp({\theta})$ is \emph{bounded}. This is the key point.

We first construct the norm appearing in (\ref{main}). Recall that the \emph{universal hypersurface associated to $B$} is given by
\begin{align}
\Sigma:= \{([F],\ [a_0:a_1:\dots :a_N])\in B\times \widehat{\cpn} | \ F(a_0,a_1,\dots, a_N)=0 \} \ .
\end{align}
Then $\Sigma$ is the base locus of the natural section
\begin{align}
\varphi \in H^0(B\times \widehat{\cpn},p_1^*\mathcal{O}_B(1)\otimes p_2^*\mathcal{O}_{\cpn}(\widehat{d})) \quad ,\ \Sigma= \{\varphi=0\} \  .
\end{align}
Let $\om$ denote the K\"ahler form on the dual $\cpn$. We consider the $(1,1)$ current $u$ on $B$ defined by the fiber integral ${p_1}_*p_2^*(\om^{N})$.
\begin{diagram}
  \Sigma &\rTo^{p_2}&\widehat{\cpn}\\
 \dTo^{p_1}\\
   B
\end{diagram}
 That is, for all $C^{\infty}$ $(b-1, b-1)$ forms $\psi$ on $B$ we have
\begin{align}\label{current}
\int_{B}u\wedge \psi=\int_{\Sigma}p_2^*(\om^N)\wedge p_1^*(\psi) \ .
\end{align}
For the following, see \cite{psc} Lemma 8.7 pg. 32 .
\begin{proposition} 
The cohomology class of the current $u$ coincides with the class of $\om_{B}$ (the Fubini-Study form). Moreover, there is a continuous function $\theta$ on $B$ such that, in the sense of currents we have
\begin{align}
u=\om_{B}+\frac{\sqrt{-1}}{2\pi}\dl\dlb \theta \ .
\end{align}
\end{proposition}

Using the identification $GX\cong G\times X$, we can exhibit the terms of the basic complex and the induced Hermitian metric $H$ in a more concrete way as follows
\begin{align*}
&K^i_{G\Delta}(m)|_{\{\sigma\}\times X}=\bigwedge^{n+1-i}T^{1,0}(\tilde{X})\otimes\mathcal{O}_{X}(m) \\
\ \\
&H^i(\sigma)=\bigwedge ^{n+1-i}(h_{\cn} \circ \sigma)|_{T^{1,0}(\tilde{X})}\otimes e^{-m\varphi_{\sigma}}h^{m}_{FS} \ .
\end{align*}
Where $h_{\cn}$ denotes the standard Hermitian form on $\cn$ .

We define a function on $G$ as follows
\begin{align}
\mathcal{I}(\sigma):= D_{K_{\Delta}^{\bull}(m)}(\mbox{Ch}_{n+1}; H^{\bull}(e), H^{\bull}(\sigma))   
\end{align}
The main point is to establish the following proposition.
\begin{proposition}\label{pluri} Let $|| \ \cdot ||:= \exp(\theta)||\ \cdot ||_{FS}$. Then the difference
\begin{align}
\mathcal{I}(\sigma)-  \log \left(\frac{{||\sigma\cdot\Delta_{X}||}^{2}}{{||\Delta_{X}||}^{2}}\right)
\end{align}
is a pluriharmonic function on $G$ .
\end{proposition}
 \noindent \emph{Proof}.  
 \begin{lemma} Let $p_i$ denote the projection onto the $ith$ factor of the incidence correspondence $I_{\Delta}$.
\begin{diagram}
  I_{\Delta} &\rTo^{p_2}&\widehat{\cpn}\\
 \dTo^{p_1}\\
  \mathbb{G}(n,\ N)
\end{diagram}
Let $\om_{\widehat{\cpn}}$ the Fubini Study K\"ahler form on $\widehat{\cpn}$. Then we have the following identity of forms on $\mathbb{G}(n,\ N)$
\begin{align}\label{pdual}
{p_1}_*(p_2^*(\om^N_{\widehat{\cpn}}))=\sum_{i=0}^{n+1}(-1)^i\mbox{\emph{Ch}}(\bigwedge^i{U}, \ h_{FS})^{\{n+1,n+1\}}\ .
\end{align}
\end{lemma}

To see this, observe that the left hand side of (\ref{pdual}) is of type $(n+1,n+1)$ and invariant under the action of the unitary group. The latter implies that it must be a polynomial in the forms $c_1(U^{\vee}),c_2(U^{\vee}),\dots, c_{n+1}(U^{\vee})$. Let $\Omega$ be any invariant form on  $\mathbb{G}(n,\ N)$ of type complimentary to ${p_1}_*p_2^*\om^N_{\widehat{\cpn}}$. Then
\begin{align}
\begin{split}
\int_{\mathbb{G}(n,\ N)}{p_1}_*(p_2^*(\om^N_{\widehat{\cpn}}))\wedge \Omega&=\int_{I_{\Delta}}p_2^*(\om^N_{\widehat{\cpn}})\wedge {p_1}^*(\Omega)\\
\ \\
&=\int_{\mathbb{G}(n,\ N)\times \widehat{\cpn}}\mbox{PD}[I_{\Delta}]\wedge p_2^*(\om^N_{\widehat{\cpn}})\wedge p_1^*(\Omega) \ .
\end{split}
\end{align}
Recall from (\ref{bslocus}) that we have $I_{\Delta}=\{s=0\}$ where $s$ is a section of $p_1^*U^{\vee}\otimes p_2^*\mathcal{O}_{\widehat{\cpn}}(+1)$. Therefore
\begin{align*}
\mbox{PD}[I_{\Delta}]&=c_{n+1}(p_1^*U^{\vee}\otimes p_2^*\mathcal{O}_{\widehat{\cpn}}(+1))\\
\ \\
&= \sum_{i=0}^{n+1}c_{1}(p_2^*\mathcal{O}_{\widehat{\cpn}}(+1))^{n+1-i}\wedge c_{i}(p_1^*U^{\vee}) \\
\ \\
&=c_{n+1}(p_1^*U^{\vee})+\sum_{i=0}^{n}c_{1}(p_2^*\mathcal{O}_{\widehat{\cpn}}(+1))^{n+1-i}\wedge c_{i}(p_1^*U^{\vee}) \ .
\end{align*}
Therefore, for \emph{all} invariant forms $\Omega$ (of complimentary type) we have
\begin{align*}
\int_{\mathbb{G}(n,\ N)}{p_1}_*(p_2^*(\om^N_{\widehat{\cpn}}))\wedge \Omega=\int_{\mathbb{G}(n,\ N)} c_{n+1}(U^{\vee})\wedge \Omega \ .
\end{align*}
Therefore,
\begin{align*}
{p_1}_*(p_2^*(\om^N_{\widehat{\cpn}}))=c_{n+1}(U^{\vee},\ h_{FS}) \ .
\end{align*}
Then the lemma follows immediately from the well known \emph{Borel-Serre identity} .
\begin{align*}
\sum_{j=0}^{k}(-1)^j\mbox{Ch}(\bigwedge^jE^{\vee})=c_{k}(E)\mbox{Td}(E)^{-1} \quad (k=rnk(E)) \ .
\end{align*}

Using the construction of the basic complex on $GX$ (see (\ref{basicmplx})) we complete the above diagram. Below $\rho_{GX}$ denotes the first component of the map $\rho_{G,f}$.

\begin{diagram}
 {\rho_{GX}}^*(I_{\Delta})  &\rTo^{\pi_2}&I_{\Delta} &\rTo^{p_2}&\widehat{\cpn}\\
\dTo ^{\pi_{1}}&&\dTo^{p_1}\\
 GX&\rTo ^{\rho_{GX}}& \mathbb{G}(n,\ N)&\\
\dTo^{\pi}\\
G
\end{diagram}
Observe that the alternating sum of the Chern Characters of the complex $(K^i_{G\Delta}(m), \ \dl_i)$ are actually \emph{independent of} $m$, and we have the identity of forms on $GX$ (where we only consider the forms of type $(n+1, n+1)$ ) .
\begin{align}
\sum_{i=0}^{n+1}(-1)^i\rho_{GX}^{*}\mbox{Ch}(\bigwedge^i U,\ h_{FS})=\sum_{i=0}^{n+1}(-1)^i\mbox{Ch}(K^i_{G\Delta}(m), \ h_{G}) \ .
\end{align}
 
 Let $\eta$ be a smooth compactly supported form on $G$ of type
 $(N^2+2N,N^2+2N)$. Then from what we have done it follows that
 \begin{align*}
 \int_G\frac{\sqrt{-1}}{2\pi}\dl\dlb \mathcal{I}\wedge \eta&=\int_{GX}\sum_{i=0}^{n+1}(-1)^{i}\mbox{Ch}_{n+1}( K^i_{G\Delta}(m), \ h^i_{G})\wedge \pi^{*}(\eta)\\
\ \\
&=\int_{GX}\rho^*_{GX}({p_1}_*(p_2^*(\om^N_{\widehat{\cpn}})))\wedge \pi^*(\eta) \\
\ \\
&= \int_{\rho^*_{GX}(I_{\Delta})}{\pi_2}^*(p_2^*(\om^N_{\widehat{\cpn}}))\wedge \pi_1^*(\pi^*(\eta)) \ .
 \end{align*}
Below $T$ denotes the evaluation map $T(\sigma):=[\sigma\cdot \Delta_{X}]$ and $\Sigma$ denotes the universal hypersurface for the family $B:=\mathbb{P}(H^0(\widehat{\cpn}, \mathcal{O}(\widehat{d})))$.
\begin{diagram}
  T^*(\Sigma)  &\rTo^{\pi_2}&\Sigma &\rTo^{p_2}&\widehat{\cpn}\\
\dTo ^{\pi_{1}}&&\dTo^{p_1}\\
 G&\rTo ^{T}& B&\\
 \end{diagram}
Let $u$ denote the positive current defined in (\ref{current}). Using the notation and commutativity in the diagram above gives that
\begin{align}
\begin{split}
\int_GT^*(u)\wedge \eta &=\int_{T^*(\Sigma)}\pi_2^*(p_2^*(\om^N_{\widehat{\cpn}}))\wedge \pi_1^*(\eta)\\
 \ \\
&=\int_{\rho^*_{GX}(I_{\Delta})}\pi_2^*(p_2^*(\om^N_{\widehat{\cpn}}))\wedge\pi^*_1\pi^*(\eta) \ .
\end{split}
\end{align}
We have used that ${T^*(\Sigma)}\cong \rho^*_{GX}(I_{\Delta}) $ (birational equivalence) . By definition we have that
\begin{align}
T^*(u)=\frac{\sqrt{-1}}{2\pi}\dl\dlb\log\left(e^{\theta\circ T}\frac{||\sigma\cdot\Delta_X||_{FS}^2}{||\Delta_X||_{FS}^2}\right) \ .
\end{align}
Therefore,
\begin{align}
\begin{split}
&\int_G\dl\dlb\left(\mathcal{I}(\sigma)-\log\left(e^{\theta\circ T}\frac{||\sigma\cdot\Delta_X||_{FS}^2}{||\Delta_X||_{FS}^2}\right)\right)\wedge \eta=0\ .\\
\end{split}
\end{align}
For all compactly supported forms $\eta$. Hence the difference is pluriharmonic.
This establishes Proposition \ref{pluri} .

 Since $G$ is simply connected there is an entire function $F$ on $G$ such that
\begin{align}
\mathcal{I}(\sigma)-\log\left(e^{\theta\circ T}\frac{||\sigma\cdot\Delta_X||_{FS}^2}{||\Delta_X||_{FS}^2}\right)
=\log(|F(\sigma)|^2) \ .
\end{align}
A standard argument shows that $F\equiv 1$. This completes the proof of the main theorem .
 \section{Smooth Plane Curves}
In this section we identify the energy for smooth curves in $\mathbb{C}P^{2}$ defined by homogeneous polynomials $F$.
 In this situation one can exploit the codimension in order to analyze and interpret the integral appearing on the left hand side of (\ref{main}). This holds for hypersurfaces of all dimensions and will be explored in a subsequent article. 
 
Precisely, we have the following expression for the ``discriminant energy" of a smooth plane curve.
\begin{align}
\begin{split}
D_{K^{\bull}}(Ch_2,H^{\bull}(e),H^{\bull}(\sigma))&=  4\deg(F)\nu_{\om}(\varphi_{\sigma}) -\deg(F)E_{1,\om}(\varphi_{\sigma})\\
\ \\
&-4\deg(\Delta_F)F_{\om}^0(\vps)\ .
\end{split}
\end{align}
 The strategy is this: instead of attempting to compute the double transgression directly \emph{we can use (\ref{cayleyenergy} ) in order to express the log of the norm of $\Delta_F$ as an integral over $X^{\vee}_F$ .} The basis for this comes from the fact that the $X$-resultant of a hypersurface $X_F$ is $F$. The integral over $X^{\vee}_F$ is then pulled back to $X_F$ by the Gauss mapping.

Let $X=X_F$ be a smooth hypersurface of dimension $n$. $F$ denotes the  irreducible defining polynomial of degree $d\geq 2$. Then the projective dual of $X_F$ is \emph{always} codimension one in 
$\widehat{\mathbb{C}P^{n+1}}$ and given by the zero set of an irreducible polynomial called the \emph{discriminant} of $F$. In this section we denote the discriminant by $\Delta_F$. The key fact is the following
\begin{align}
\rho_{F}(X_F)=\widehat{X_F}
\end{align}
 $\rho_F$ denotes the \emph{Gauss map} of $F$ and is given explicitly by the following formula.
\begin{align*}
\rho_F(p)=[\frac{\dl F}{\dl z_0}(p): \frac{\dl F}{\dl z_1}(p):\frac{\dl F}{\dl z_2}(p):\dots:\frac{\dl F}{\dl z_{n+1}}(p)]\  ; \ p\in X_F
\end{align*}
 It is well known that $\rho_F$ is a \emph{birational isomorphism}. In fact $\rho_F:X_F\rightarrow X_F^{\vee}$ is a \emph{resolution of singularities} .
We denote coordinates on the dual projective space by $[a]=[a_0:a_1:\dots:a_{n+1}]$. Recall that the dual action of $GL(n+2,\mathbb{C})$ on $[a]$ is given by
\begin{align}\label{dualaction}
\sigma\cdot a:= (\sigma^{-1})^t a
\end{align}
$a$ is viewed as a column vector and the right hand side of (\ref{dualaction}) is just matrix multiplication.
The corresponding \emph{dual Bergman potential} is given by the formula
\begin{align*}
\widehat{\varphi_{\sigma}}([a])=\log\frac{|\sigma\cdot a|^2}{|a|^2} \ .
\end{align*}

Recall the well known fact (see \cite{kenhyp}) .
\begin{proposition}\label{riccicurv}
\begin{align} 
\begin{split}
&\mbox{Ric}(\om |_{X_F})=(n+2-d)\om-\frac{\sqrt{-1}}{2\pi}\dl\dlb\psi_{F}\\
\ \\
&\psi_F(z):= \log\left(\frac{\sum_{i=0}^{n+1}|\frac{\dl F}{\dl z_i}|^2}{||z||^{2(d-1)}}\right)\ .
\end{split}
\end{align}
\end{proposition}
\begin{proof}

By definition of the Gauss map we have
 \begin{align*}
 \rho_{F}^{*}(\mathcal{U})\cong T^{1,0}(\tilde{X}_F) \ .
 \end{align*}
Where $\mathcal{U}$ denotes the universal bundle over $\mathbb{G}(n\ , \mathbb{C}P^{n+1})$.
Observe that via the natural isomorphism $\iota$
\begin{align*}
\mathbb{G}(n\ , \mathbb{C}P^{n+1})\overset{\iota}{\cong} \widehat{\mathbb{C}P^{n+1}}
\end{align*}
we have the identification
\begin{align*}
\bigwedge^{n+1}\mathcal{U}\cong \iota^*\mathcal{O}_{\cpnd}(-1) \ .
\end{align*}
Therefore,
\begin{align}
c_1\left(\bigwedge^{n+1} T^{1,0}(\tilde{X})\right) = - \rho_{F}^{*}\omd\ .
 \end{align}
Recall that the dual Fubini study form is given by
\begin{align*}
\omd=\frac{\sqrt{-1}}{2\pi}\dl\dlb \log (|a_0|^2+|a_1|^2+\dots |a_{n+1}|^2) \ .
\end{align*}
 
Therefore we have that
\begin{align}\label{pllbck}
-c_1\left(\bigwedge^{n+1} T^{1,0}(\tilde{X})\right)=(d-1)\om+\frac{\sqrt{-1}}{2\pi}\dl\dlb\log\left(\frac{\sum_{i=0}^{n+1}|\frac{\dl F}{\dl z_i}|^2}{||z||^{2(d-1)}}\right)\ .
\end{align}

Recall the exact sequence
\begin{align*}
 0\rightarrow \mathcal{O}_X(-1)\overset{\iota}{\rightarrow}&T^{1,0}(\tilde{X})\ \overset{\pi}{\rightarrow} T^{1,0}(X)\otimes \mathcal{O}_X(-1)\rightarrow 0 \ .
\end{align*}
All the terms of this sequence are equipped with natural Hermitian metrics induced from $T^{1,0}(\tilde{X})$.
Recall the general \emph{curvature decomposition}
\begin{align*}
 F_\mathcal{E}= \begin{pmatrix}
 F_\mathcal{S}-\beta^*\wedge \beta & D^{1,0}_{Hom(\mathcal{Q},\mathcal{S})}\beta^*\\
\dlb \beta& F_\mathcal{Q}-\beta\wedge \beta^*
\end{pmatrix}
\end{align*}
associated to any short exact sequence of holomorphic Hermitian vector bundles.
\begin{align*}
0\rightarrow \mathcal{S}\rightarrow \mathcal{E}\rightarrow \mathcal{Q}\rightarrow 0 \ .
\end{align*} 
$\beta^*\in C^{\infty}(Hom(\mathcal{Q},\mathcal{S})\otimes \Omega^{0,1}_X)$ denotes the \emph{second fundamental form} of the inclusion $0\rightarrow \mathcal{S}\rightarrow \mathcal{E}$ .
In our case we have that
\begin{align*}
Tr_{\mathcal{S}}(\beta^*\wedge \beta)+Tr_{\mathcal{Q}}(\beta\wedge \beta^*) =0 \ .
\end{align*}
This implies that
\begin{align*}
c_1\left(\bigwedge^{n+1} T^{1,0}(\tilde{X})\right)= -(n+1)\om +\mbox{Ric}(\om|_{X_F}) \ .
\end{align*}
Putting this together with (\ref{pllbck}) proves the proposition.
\end{proof}
 
\begin{proposition}\label{ddbar}
For all $\sigma \in GL(n+2,\mathbb{C})$ we have
\begin{align*}
\widehat{\varphi_{\sigma}}\circ \rho_F(z)=(n+1)\vps(z)+\log\left(\frac{\om^n_{\sigma}}{\om^n}\right)-
\log(|\mbox{\emph{det}}(\sigma)|^2) \ .
\end{align*}
\end{proposition}
We begin by establishing the following identity. In the statement we have defined $F^{\sigma}:= \sigma\cdot F$ .
\begin{lemma} There is a function $C: GL(n+2,\mathbb{C})\ni\sigma \ra C(\sigma)\in\mathbb{R}$ such that 
\begin{align}
(n+2-d)\varphi_{\sigma}+\log\left(\frac{\om^n_{\sigma}}{\om^n}\right)=\psi_{F^{\sigma}}(\sigma z)-\psi_F(z) + C(\sigma) \ .
\end{align}
\end{lemma}
\begin{proof}
The argument consists in computing $Ric(\oms|_{X_F})$ in two different ways. 
Obviously $\sigma X_F=X_{F^{\sigma}}$. Therefore,
\begin{align}\label{ricsigma}
Ric(\om|_{\sigma X_F})=Ric(\om|_{X_{F^{\sigma}}})  \ .
\end{align}
Since $\sigma^{*}(\om|_{\sigma X_F})=\om_{\sigma}|_{X_F}$ we have that 
\begin{align*}
Ric(\oms|_{X_F})= Ric(\sigma^*\om|_{X_{F^{\sigma}}})=\sigma^*Ric(\om|_{X_{F^{\sigma}}})  \ .
\end{align*}
Therefore by proposition (\ref{riccicurv}) we see that 
\begin{align}\label{restriction}
Ric(\oms|_{X_F})=(n+2-d)\oms|_{X_F}-\frac{\sqrt{-1}}{2\pi}\dl\dlb \psi_{F^{\sigma}}\circ\sigma \ .
\end{align}
On the other hand the definition of the Ricci form gives at once that 
\begin{align}\label{defricci}
\begin{split}
Ric(\oms|_{X_F})&=-\frac{\sqrt{-1}}{2\pi}\dl\dlb \log\left(\frac{\om^n_{\sigma}}{\om^n}\right)+Ric(\om)\\
\ \\
&=-\frac{\sqrt{-1}}{2\pi}\dl\dlb \log\left(\frac{\om^n_{\sigma}}{\om^n}\right)-\frac{\sqrt{-1}}{2\pi}\dl\dlb \psi_F+(n+2-d)\om \ .
\end{split}
\end{align}
Combining (\ref{restriction}) and (\ref{defricci}) we deduce the following $\dl\dlb$ equation
\begin{align*}
\dl\dlb \left((n+2-d)\varphi_{\sigma}+\log\left(\frac{\om^n_{\sigma}}{\om^n}\right)-\psi_{F^{\sigma}}(\sigma z)+\psi_F(z)\right)=0 \ .
\end{align*}
Therefore there is a constant $C(\sigma)$ such that
\begin{align*}
(n+2-d)\varphi_{\sigma}+\log\left(\frac{\om^n_{\sigma}}{\om^n}\right)=\psi_{F^{\sigma}}(\sigma z)-\psi_F(z) +C(\sigma) \ .
\end{align*}
\end{proof}
By definition of the dual potential we have 
\begin{align*}
\widehat{\vp}_{\tau}\circ \rho_F= \psi_{F^{\sigma}}(\sigma z)-\psi_F(z) +(d-1)\vps \ .
\end{align*}
Combining this with proposition (\ref{ddbar}) gives
\begin{align}\label{dualpotlog}
\widehat{\vp}_{\sigma}\circ \rho_F=(n+1)\vps+\log\left(\frac{\om^n_{\sigma}}{\om^n}\right)-C(\sigma) \ .
\end{align}
\begin{claim}
For all $\sigma$ and $\tau$ in $GL(n+2, \mathbb{C})$ we have
\begin{align*}
C(\sigma\tau)=C(\sigma)+C(\tau)\ . \\
\end{align*}
Therefore $C$ is a homomorphism from $GL(n+2, \mathbb{C})$ into the additive group $\mathbb{R}$.
\end{claim}
 
\begin{proof} Observe that
\begin{align*}
\begin{split}
\vps \circ \tau=\log\left(\frac{|\sigma \tau z|^2}{|\tau z|^2}\right)=\log\left(\frac{|\sigma \tau z|^2}{|z|^2}\right)
-\log\left(\frac{|\tau z|^2}{| z|^2}\right)=\varphi_{\sigma \tau}(z)-\varphi_{\tau}(z)
\end{split}
\end{align*}
Therefore,
\begin{align}\label{potident}
  \vps \circ \tau+ \varphi_{\tau}(z)= \varphi_{\sigma \tau}(z) \ .
\end{align}
Which in turn implies that
\begin{align*}
\tau^*(\oms)=\om_{\sigma\tau} \ .
\end{align*}
Therefore,
\begin{align}\label{logident}
\log\left(\frac{\om^n_{\sigma}}{\om^n}\right)\circ \tau=\log\left(\frac{\om^n_{\sigma\tau}}{\om^n}\right)-\log\left(\frac{\om^n_{\tau}}{\om^n}\right) \ .
\end{align}
Similarly we have 
\begin{align}\label{dualberg}
\widehat{\vp}_{\sigma\tau}=\widehat{\vp}_{\tau}+\widehat{\vps}\circ \tau \ .
\end{align}
An application of the chain rule gives
\begin{align}\label{chain}
\tau\cdot \rho_F= \rho_{F^{\tau}}\circ\tau \ .
\end{align}
(\ref{dualberg}) and (\ref{chain}) imply that
\begin{align*}
\begin{split}
\widehat{\vp}_{\sigma\tau}\circ\rho_F&=\widehat{\vp}_{\tau}\circ \rho_F+\widehat{\vps}\circ \tau\cdot\rho_F\\
\ \\
&=\widehat{\vp}_{\tau}\circ \rho_F+\widehat{\vps}\circ \rho_{F^{\tau}}\circ\tau \ .
\end{split}
\end{align*}
Combining this and (\ref{dualpotlog}) shows that
\begin{align}\label{bigidentity}
\begin{split}
&(n+1)\vp_{\sigma\tau}+\log\left(\frac{\om^n_{\sigma\tau}}{\om^n}\right)-C(\sigma\tau)=\\
\ \\
&(n+1)(\vp_{\tau}+\vps\circ\tau)+\log\left(\frac{\om^n_{\sigma}}{\om^n}\right)\circ\tau+\log\left(\frac{\om^n_{\tau}}{\om^n}\right)-(C(\sigma)+C(\tau)) \ .
\end{split}
\end{align}
Apply (\ref{logident}) and (\ref{potident}) in order to finish the proof of the claim. \end{proof}

It is easy to see that on all matrices of the form $t\mathbb{I}_{n+2}$ (where $t\in \mathbb{C}^*$) we have
\begin{align}
\exp C(t\mathbb{I}_{n+2})=|t|^{2(n+2)}=|\mbox{det}(t\mathbb{I}_{n+2})|^2 \ .
\end{align}
The claim shows that $C$ is a \emph{class function}. Therefore, 
\begin{align*}
C((t,1,\dots,1))=C((1,t,1,\dots 1)=\dots = C((1,1,\dots,t)) \ .
\end{align*}
Therefore
\begin{align*}
\left(\exp C(t,1,\dots,1)\right)^{n+2}=|t|^{2(n+2)}\ .
\end{align*}
Therefore on all \emph{diagonal} matrices we have
\begin{align*}
\exp C(t_1,t_2,\dots,t_{n+2})=|t_1|^2|t_2|^2\dots |t_{n+2}|^2 \ .
\end{align*}
Since diagonalizable matrices are \emph{dense} in $GL(n+2, \mathbb{C})$ and $C$ depends continuously on $\sigma$ we conclude that $C(\sigma)=\log|\mbox{det}(\sigma)|^2$ .  
 
 Now we come to the main result in this section which expresses the dual variety and the chow point in terms of the Liouville energy and the Mabuchi energy restricted to the Bergman metrics. As mentioned in the beginning of this section
we do not compute $D_{K^{\bull}}$ directly. Instead we obtain an energy expression for $\log||\sigma\cdot\Delta_F||$ through an application of (\ref{cayleyenergy}). 
\begin{proposition}\label{insteadweobtain}
Let $X_F$ be a smooth hypersurface in $\mathbb{C}P^2$. Let $\Delta_F$ denote the discriminant of $F$. Then
\begin{align} 
\begin{split}
 4\deg(F)\nu_{\om}(\varphi_{\sigma})&= \log \left(\frac{{||\sigma\cdot\Delta_{F}||}^{2}}{{||\Delta_{F}||}^{2}}\right)- 2\frac{\deg(\Delta_F)}{\deg(F)}\log \left(\frac{{||\sigma\cdot{F}||}^{2}}{||{F}||^2}\right)\\
 \ \\
 &+\deg(F)E_{1,\om}(\varphi_{\sigma}) \ . 
\end{split}
\end{align}
\end{proposition}
\begin{proof}
We have the following immediate application of (\ref{cayleyenergy}) .
\begin{align}  \log \left(\frac{{||\sigma\cdot\Delta_{F}||}^{2}}{{||\Delta_{F}||}^{2}}\right)=-\int_{X_F^{\vee}}\dl\widehat{\vp}_{\sigma}\wedge \dlb\widehat{\vp}_{\sigma}+2\int_{X_F^{\vee}}\widehat{\vp}_{\sigma}\widehat{\om} \end{align}
Next use that ${X_F^{\vee}}=\rho_F(X_F)$ is the \emph{birational image} of $X_F$ under the Gauss map $\rho_F$ .
Pulling everything back to $X_F$ gives the identity
\begin{align} \label{trick}
\log \left(\frac{{||\sigma\cdot\Delta_{F}||}^{2}}{{||\Delta_{F}||}^{2}}\right)= -\int_{X_F}\dl(\widehat{\vp}_{\sigma}\circ\rho_F)\wedge \dlb(\widehat{\vp}_{\sigma}\circ\rho_F)+2\int_{X_F^{\vee}}(\widehat{\vp}_{\sigma}\circ\rho_F)\rho^*_F(\widehat{\om})\ .
\end{align}
 Since $n=1$ and $\sigma \in SL(3,\mathbb{C})$ proposition (\ref{ddbar}) implies at once that
\begin{align}\label{above1}
\widehat{\vp}_{\sigma}\circ\rho_F=2\vps + \log\left( \frac{\oms}{\om}\right) \ .
\end{align}
Substituting (\ref{above1}) into the first integral on the right hand side of (\ref{trick}) gives
\begin{align*} 
\begin{split}
\int_{X_F}\dl(\widehat{\vp}_{\sigma}\circ\rho_F)\wedge \dlb(\widehat{\vp}_{\sigma}\circ\rho_F)&=4\int_{X_F}\dl\vps\wedge\dlb\vps-4\int_{X_F}\log\left( \frac{\oms}{\om}\right)\oms \\ 
\ \\
&+4\int_{X_F}\log\left( \frac{\oms}{\om}\right)\om +\int_{X_F}|\nabla \log\left( \frac{\oms}{\om}\right)|^2\om \ .
\end{split}
\end{align*}
By definition of the K-Energy we have (up to a bounded term which we ignore)
\begin{align}
4\int_{X_F}\log\left( \frac{\oms}{\om}\right)\oms=4\deg(F)\nu_{\om}(\vps)+2(3-\deg(F))\int_{X_F}\dl\vps\wedge\dlb\vps \ .
\end{align}
Therefore we have the following identity ,
\begin{align}
\begin{split}
-\int_{X_F}\dl(\widehat{\vp}_{\sigma}\circ\rho_F)\wedge \dlb(\widehat{\vp}_{\sigma}\circ\rho_F)&= -2(\deg(F)-1)\int_{X_F}\dl\vps\wedge\dlb\vps \\
\ \\
&+ 4\deg(F)\nu_{\om}(\vps)+  4\int_{X_F}\psi_F\oms\\
\ \\
& -4\int_{X_F}\log\left( \frac{\oms}{\om}\right)\om
 -\int_{X_F}|\nabla \log\left( \frac{\oms}{\om}\right)|^2\om \ .
\end{split}
\end{align}
 By definition of the Gauss map we have
\begin{align*} 
\rho^*_F(\widehat{\om})=(\deg(F)-1)\om+\dl\dlb\psi_F \ . 
\end{align*}
Combining this with proposition \ref{ddbar} gives the following expression for the mean
\begin{align}
\begin{split}
2\int_{X^{\vee}_F}\widehat{\vp}_{\sigma} \widehat{\om}&= 4(\deg(F)-1) \int_{X_F}\vps\om + 4\int_{X_F}\vps \dl\dlb\psi_F\\
\ \\
&+2(\deg(F)-1) \int_{X_F}\log\left( \frac{\oms}{\om}\right)\om \\
\ \\
&+2\int_{X_F}\log\left( \frac{\oms}{\om}\right)\dl\dlb\psi_F \ .
\end{split}
\end{align}
 Putting all of this together gives the following identity.
\begin{align}\label{thistogether}
\begin{split}
  \log \left(\frac{{||\sigma\cdot\Delta_{F}||}^{2}}{{||\Delta_{F}||}^{2}}\right)&= 4\deg(F)\nu_{\om}(\vps)-4\deg(F)(\deg(F)-1)F_{\om}^0(\vps)\\
 \ \\
 & +2(\deg(F)-3)\int_{X_F}\log\left( \frac{\oms}{\om}\right)\om -\int_{X_F}|\nabla \log\left( \frac{\oms}{\om}\right)|^2\om\\
\ \\
&  + 2\int_{X_F}\log\left( \frac{\oms}{\om}\right)\dl\dlb\psi_F \ .
  \end{split}
  \end{align}
 Next apply proposition \ref{riccicurv} to (\ref{thistogether}) in order to get
\begin{align}\label{inorder}
 \begin{split}
 \log \left(\frac{{||\sigma\cdot\Delta_{F}||}^{2}}{{||\Delta_{F}||}^{2}}\right)= &4\deg(F)\nu_{\om}(\vps)-4\deg(F)(\deg(F)-1)F_{\om}^0(\vps)\\
& -2\int_{X_F}\log\left( \frac{\oms}{\om}\right)Ric(\om)
 -\int_{X_F}|\nabla \log\left( \frac{\oms}{\om}\right)|^2\om \ .
 \end{split}
 \end{align}
 Integration by parts yields
\[-\int_{X_F}|\nabla \log\left( \frac{\oms}{\om}\right)|^2\om=\int_{X_F}\log\left( \frac{\oms}{\om}\right)(Ric(\om)-Ric(\oms))\ . \] 
 Therefore
 \begin{align*}
 2\int_{X_F}\log\left( \frac{\oms}{\om}\right)Ric(\om)
 +\int_{X_F}|\nabla \log\left( \frac{\oms}{\om}\right)|^2\om =\deg(F)E_{1,\om}(\vps) \ .
 \end{align*}
 Applying this to (\ref{inorder}) gives
 \begin{align}\label{applyingthis}
 \begin{split}
 \log \left(\frac{{||\sigma\cdot\Delta_{F}||}^{2}}{{||\Delta_{F}||}^{2}}\right)&= 4\deg(F)\nu_{\om}(\vps)-4\deg(F)(\deg(F)-1)F_{\om}^0(\vps)\\
 &- \deg(F)E_{1,\om}(\vps) \ .
 \end{split}
 \end{align}
 It is well known that
\begin{align}\label{degree}
{\deg}(\Delta_F)={\deg}(F)({\deg}(F)-1) \ . 
\end{align}
Appealing to (\ref{cayleyenergy}) once more yields
\begin{align*}
-2\deg(F)F_{\om}^0(\vps)=\log\left(\frac{||\sigma\cdot F||^2}{||F||^2}\right) \ .
\end{align*}
Direct substitution of this into (\ref{applyingthis}) together with (\ref{degree}) completes the proof of proposition \ref{insteadweobtain}.
\end{proof}
\section{Further Remarks}
 
 In this subsection we discuss some of the consequences of our working conjecture. Let $E$ be a finite dimensional complex rational representation of a  complex torus $H:=(\mathbb{C}^{*})^{N}$. As usual $\chi (H)$ denotes the character group of $H$
\begin{align*}
\chi  (t_{1},t_{2},\dots,,t_{N})=t_{1}^{m_{1}}t_{2}^{m_{2}}\dots t_{N}^{m_{N}} \ , \qquad m_{i} \in \mathbb{Z}\ .
\end{align*}
$\chi (H)$ is lattice of full rank in the finite dimensional real vector space $\chi (H)\otimes_{\mathbb{Z}} \mathbb{R}$.\newline
$E$ decomposes under the $H$ representation into \emph{weight spaces} $E_{\chi}$
\begin{align*}
E= \bigoplus_{\chi \in \chi (H)}E_{\chi} \qquad t\in H \ \mbox{acts on $E_{\chi}$ by $\chi(t)$ } \ .
\end{align*}
Let $v\in E$ be a nonzero vector in $E$ then $v$ decomposes into weight vectors
\begin{align}\label{supp}
v=\sum _{\chi \in\ \mbox{supp}(v)}v_{\chi} \ .
\end{align}
supp($v$) denotes the \emph{support} of $v$.  supp($v$) consists of all $\chi \in \chi(H)$ such that $v_{\chi}\neq 0$ (the projection of $v$ into $E_{\chi}$).
 
 A \emph{one parameter subgroup} of $G$ is an algebraic\footnote{``algebraic'' means that the matrix coefficients $\lambda(t)_{i,j}\in \mathbb{C}[t,t^{-1}]$.} homomorphism
\begin{align*}
\lambda:\mathbb{C}^{*}\rightarrow  G\ .
\end{align*}
 Any such $\lambda(t)$ can be diagonalised. That is, we may assume that $\lambda(t)$ takes values in the standard maximal torus $H\cong (\mathbb{C}^*)^N$ of $G$.   
\begin{align*}
\lambda(t)=\begin{pmatrix}t^{m_{0}}&\dots&\dots& 0\\
                             0&t^{m_{1}}&\dots& 0\\
                             0&\dots&\dots& t^{m_{N}}
                             \end{pmatrix}\ .\end{align*}
The exponents $m_{i}$ satisfy
\begin{align*}
 \quad \sum_{0\leq i \leq N}m_{i}=0.
\end{align*}
The space of one parameter subgroups will be denoted by $\Gamma (H)$ .
Then, following the considerations of section 3 we have
 \begin{align*}
      \lambda(t)^{*}\omega_{FS}|_{X}=\omega +\frac{\sqrt{-1}}{2\pi}\dl\dlb\mbox{log}\left(\sum_{0\leq j\leq N}|t|^{2m_{j}}||S_{j}||^{2}(z)\right)\ .
\end{align*}
 Below we shall abuse notation somewhat and define
 \begin{align*}
 \varphi_{ \lambda(t)}:= \mbox{log}\left(\sum_{0\leq j\leq N}|t|^{2m_{j}}||S_{j}||^{2}(z)\right)\ .
 \end{align*}

 Recall that the dual to the space of characters is the space of {one parameter subgroups} .
  The duality is given as follows
\begin{align*}
\chi (\lambda(t))= t ^{<\chi,\lambda>}:= t^{m_{1}a_{1}+\dots + m_{N}a_{N}} \ .
\end{align*}
In other words there is an isomorphism
\begin{align*}
\Gamma(H)\otimes_{\mathbb{Z}}\mathbb{R}\cong (\chi(H)\otimes_{\mathbb{Z}}\mathbb{R})^{\vee}\ .
\end{align*}

 Let $P(v)$ denote the convex hull of all $\chi \in \mbox{supp}(v)$. Then $P(v)$ is a  compact convex integral polytope (the \emph{weight polytope}) inside $\chi(H)\otimes_{\mathbb{Z}}\mathbb{R}$.
 The integral linear functional corresponding to $\lambda$ will be denoted by $l_{\lambda}$.
\begin{definition}
The weight  $w_{\lambda}(v)$  of $\lambda$ on $v\in E$ is the integer
\begin{align*}
w_{\lambda}(v):= \mbox{\emph{Min}}_{ \{ x\in P(v)\}}\ l_{\lambda}(x)= \mbox{\emph{Min}} \{ <\chi,\lambda>| \chi \in \mbox{\emph{supp}}(v)\}
\end{align*}
\end{definition}
It is clear that  $w_{\lambda}(v)$ is the unique integer such that
\begin{align*}
\lim_{t\rightarrow 0}t^{-w_{\lambda}(v)}\lambda(t)v \quad \mbox{\emph{exists and is \textbf{not} zero}}.
\end{align*}
Let $|| \ ||$ denote any norm on $E$. Then we have that
\begin{align}\label{kempfness}
\lim_{t \rightarrow 0}\log (||\lambda(t)v||^{2})=w_{\lambda}(v)\log(|t|^2)+O(1)\ .
\end{align}
Our working hypothesis implies at once that there is an expansion as $|t|\rightarrow 0$
\begin{align}\label{asymptotic}
\nu_{\om}(\vplt)=\left(\kappa_1w_{\lambda}(v_1)-\kappa_2w_{\lambda}(v_2)\right)\log(|t|^2)+O(1) \ .
\end{align}

\begin{definition} (Tian \cite{psc})
{$\nu_{\omega}$ is \textbf{proper} if there exists a strictly increasing function $f:\mathbb{R}_{+}\ra\mathbb{R}_{+}$ (where $\lim_{T\ra \infty}f(T)=\infty$) such that $\nu_{\omega}(\varphi)\geq f(J_{\omega}(\varphi))$ for all $\varphi\in  \mathcal{H}_\om$}.
 \end{definition}

We have the following corollary of our working conjecture.  
Below, $\eta (X)$ denotes the space of holomorphic vector fields on $X$.\\
 \ \\
 \textbf{Corollary of Conjecture 1 .}\emph{ Assume that $\eta(X)=\{ 0 \}$. Then
 the Mabuchi energy is proper along all degenerations $\lambda\in \Gamma(H)$ if and only if there is a positive constant $C=C(\om)$ such that}
 \begin{align} \label{mult}
 \kappa_1w_{\lambda}(v_1)-\kappa_2w_{\lambda}(v_2) + \frac{C}{{\deg}(X)(n+1)}e(\lambda; \ X)\leq 0\ .   
 \end{align}
 \emph{$ e(\lambda; \ X)$ denotes the \textbf{multiplicity} (see definition 2.2 pg. 55 of \cite{sopv}) of $X$ with respect to $\lambda$. Moreover, in this case, the scaled weight polytope of $v_2$ strictly dominates the weight polytope of $v_1$ . The scaling factor being} $\frac{\kappa_2}{\kappa_1}$ .
\begin{remark}{ (\ref{mult}) was inspired by a lecture of David Calderbank in March 2008 at the DeGiorgi institute . }
\end{remark}
\begin{center}{\textbf{Acknowledgments}}\end{center}
The author is in debt to Eckart Viehweg who suggested considering the dual variety. Jeff Viaclovsky provided   criticism of several early drafts of this paper, his many suggestions improved the quality of the exposition. This work was supported by a NSF DMS grant 0505059 .
\bibliography{ref}
\end{document}